\documentclass[english,12pt]{article}
\usepackage[T1]{fontenc}
\usepackage[utf8]{inputenc}

\usepackage{fullpage}
\usepackage{amsmath,amsthm}
 
\newtheorem{thm}{Theorem}
\numberwithin{thm}{section}
\newtheorem{lem}{Lemma}
\numberwithin{lem}{section} 

\theoremstyle{definition}

\numberwithin{defi}{section}

\numberwithin{exper}{section}
\theoremstyle{exper}

\theoremstyle{definition}
\newtheorem{assumption}{Assumption}
\numberwithin{assumption}{section}

\theoremstyle{definition}

\numberwithin{exam}{section}

\numberwithin{hequa}{section} 

\newtheorem{rem}{Remark}

\usepackage{amssymb}
\usepackage{graphicx} 
\usepackage{subfigure}
\usepackage{color}
\usepackage{epsfig}
\usepackage{epstopdf}
\usepackage{hyperref}
\usepackage{algorithm,algcompatible}
\usepackage{authblk}
\usepackage{comment}
\usepackage{xcolor}

\numberwithin{equation}{section}
\usepackage{booktabs}

\graphicspath{ {./figures/} }

\title{Essentially Order-Preserving Fourth-Order Exponential Runge--Kutta Methods for Non-Commuting Operators}
	
\author[1]{Thi Tam Dang}
\author[2]{Trung Hau Hoang\thanks{\href{mailto:trung-hau.hoang@matfyz.cuni.cz}{Corresponding author: trung-hau.hoang@matfyz.cuni.cz (Trung Hau Hoang)\\
			Email address Thi Tam Dang:
			 tam.dang@helsinki.fi
			 }}}
\affil[1]{Department of Mathematics and Statistics, University of Helsinki, Finland} 
\affil[2]{Faculty of Mathematics and Physics, Charles University, Czech Republic}

\date{}
\begin{document}

\maketitle

\begin{abstract}
This paper investigates the convergence of explicit exponential
Runge--Kutta methods for the linear evolution equation
\(u'(t)+Au(t)=Bu(t)\), where \(-A\) generates an analytic semigroup,
\(B\) is unbounded and relatively bounded with respect to \(A\), and
\(A\) and \(B\) do not commute. It is known from our previous work on
order reduction of exponential Runge--Kutta methods for non-commuting
operators that this non-commutativity can cause order reduction even for
linear problems. In a further work on third-order methods, we showed that
this reduction can be avoided by adding a stage and imposing higher
regularity on the initial data. Here we treat the fourth-order case. We
identify conditions on the coefficients of explicit fourth-order
exponential Runge--Kutta methods under which the full order four is
retained for this class of problems, and we state the convergence rate as
a function of the regularity of \(u_0\). The analysis shows precisely
which defect terms, caused by the non-commutativity of \(A\) and \(B\),
are eliminated by these conditions. The classical methods of Krogstad and
of Strehmel--Weiner do not satisfy them and converge only with order
\(11/4\) under the corresponding regularity assumptions, whereas the
six-stage method expRK4s6, constructed by Luan (2021), satisfies them and
attains order four. Numerical experiments for an advection--diffusion problem involving non-commuting operators confirm the theoretical convergence rates.

\end{abstract}



\section{Introduction}\label{sec:intro}

Exponential integrators are well suited to stiff evolution equations, since the stiff linear part is treated exactly and only the remainder is integrated explicitly; see \cite{HO2010} for a survey. For semilinear problems with a bounded nonlinearity, the convergence theory of exponential Runge--Kutta methods is well developed \cite{HO2005}. The situation changes when the explicitly treated term is itself unbounded and does not commute with the stiff operator, as in advection--diffusion or convection--reaction--diffusion models. We study the linear model problem
\begin{equation}\label{eq:model}
	u'(t) + Au(t) = Bu(t), \qquad u(0)=u_0,
\end{equation}
in a Banach space $X$, where $-A$ generates an analytic semigroup, $B$ is relatively bounded with respect to $A$, and in general $AB\neq BA$. A typical example is a second-order elliptic operator $A$ combined with a first-order operator $B$. Problem \eqref{eq:model} is a linear test for the nonlinear parabolic equation $\partial_t u-\Delta u=f(\nabla u,u)$; extending the analysis to that case requires additional nonlinear estimates and is not attempted here.

In \cite{Hau2025} we showed that for \eqref{eq:model} the non-commutativity of $A$ and $B$ causes order reduction: explicit exponential Runge--Kutta methods that satisfy the stiff order conditions only in a weak form (that is, with $A=0$) lose part of their classical order. For third-order methods the observed rate drops to about $5/2$, and for the fourth-order methods of Krogstad and of Strehmel--Weiner the convergence order we can prove is $11/4$ (see \cite{Hau2025,hau2026order}). In \cite{HauThird} we showed that, at third order, this reduction can essentially be removed by a method that satisfies the stiff order conditions in their strong form, under stronger regularity of $u_0$.

Here we treat fourth-order methods. The main contributions are:
\begin{enumerate}
	\item We identify order conditions on the coefficients of explicit fourth-order exponential Runge--Kutta methods which guarantee essentially order four for \eqref{eq:model} with non-commuting unbounded $A$ and $B$. The result is stated for a class of methods and is not tied to one scheme.
	\item The error recursion shows which defect terms caused by $AB\neq BA$ vanish.
	\item Under  $u_0\in D(A-B)^{3}$, we prove $\|u_n-u(t_n)\|
	\leq
	C t_n^{-1+2\zeta}
	\tau^{4-2\zeta},
	 \text{uniformly for }1\leq n\leq N,$ with $C$ independent of $\tau$ and $n$.

	\item The six-stage method expRK4s6 of Luan \cite{Luan2021} satisfies the conditions and therefore attains order nearly four, in contrast with the order $11/4$ obtained for the Krogstad and Strehmel--Weiner methods.
\end{enumerate}
For the convergence analysis, we impose an additional structural
condition on \(A\), following the approach adopted in the
aforementioned works. More precisely, we assume that \(A\) admits
a bounded \(H^\infty\)-functional calculus on a suitable sector.
This condition allows the operator norm of \(f(A)\) to be
bounded in terms of the supremum norm of the scalar function
\(f\) on that sector. Numerical experiments for an
advection--diffusion problem with variable coefficients are also provided.

The paper is organized as follows. Section~\ref{sec:framework} presents the analytical framework. Section~\ref{recallEXP} recalls exponential Runge--Kutta methods and the associated stiff order conditions. Section~\ref{erroranalysis} develops the error analysis, derives the error recursion, and establishes the main convergence result. Numerical experiments are presented in Section~\ref{numericalchap2}, followed by concluding remarks in Section~\ref{concluchap2}. 

\section{Analytical framework}\label{sec:framework}

We work in the framework of \cite{Hau2025,HauThird} and recall only what is needed here. Let $X$ be a Banach space with norm $\|\cdot\|$, and let $\mathcal{D}(A)$ denote the domain of $A$ in $X$. For the semigroup theory we refer to \cite{henry1981geometric,pazy1983semigroups}.

\begin{assumption}
	\label{ass1}
		The operator $A:\mathcal{D}(A)\to X$ is sectorial and invertible, i.e.\ $0\in\rho(A)$.
	Assume that there exist \(\theta\in(0,\pi)\) and \(C_A>0\) such that \(A\)
	admits a bounded \(H^\infty(\Omega_\theta)\)-functional calculus, where
	\[
	\Omega_\theta:=\{z\in\mathbb C\setminus\{0\}:|\arg z|<\theta\}.
	\]
	That is, the map \(f\mapsto f(A)\) extends to a bounded linear map from
	\(H^\infty(\Omega_\theta)\) into \(\mathcal L(X)\), and
	\[
	\|f(A)\|\le C_A\|f\|_{H^\infty(\Omega_\theta)}
	\qquad\text{for all }f\in H^\infty(\Omega_\theta).
	\]
\end{assumption}

Under Assumption~\ref{ass1}, $-A$ is the infinitesimal generator of an analytic semigroup $\{\mathrm{e}^{-tA}\}_{t\ge 0}$, and the fractional powers $A^{\gamma}$, $\gamma\ge 0$, are well defined.  The bounded \(H^\infty\)-functional calculus allows
us to transfer scalar bounds for holomorphic functions to operator-norm bounds
for their functional-calculus realizations.

The following result shows that a first-order cancellation of the quadrature
defect at the origin, together with uniform boundedness away from the origin,
yields a power-type decay estimate for the associated operator. This estimate
is useful in the convergence analysis because it quantifies the contribution
of the quadrature defect: terms of the form
\(\chi(-\tau A)A^{-\mu}u\) decay like \(\tau^\mu\) as \(\tau\to0\).
Consequently, if \(u\in\mathcal D(A^\mu)\), one obtains
\[
\|\chi(-\tau A)u\|
\le
C\tau^\mu\|A^\mu u\|,
\]
so regularity of the data measured by \(\|A^\mu u\|\) translates into an
explicit convergence rate.

\begin{thm}
	\label{thm:quadrature-defect}
	Let \(A\) satisfy Assumption~\ref{ass1}.
	Let \(\chi\) be holomorphic on a sector containing \(-\tau\Omega_\theta\),
	and assume that
	\[
	\chi(0)=0.
	\]
	Suppose further that, for some \(C_\chi>0\),
	\[
	\sup_{w\in-\tau\Omega_\theta}|w^{-\mu}\chi(w)|\le C_\chi
	\qquad\text{for every }0\le\mu\le1.
	\]
	Then, for every \(0\le\mu\le1\),
	\[
	\|\chi(-\tau A)A^{-\mu}\|\le C_A C_\chi\,\tau^\mu.
	\]
	In particular, if \(C_\chi\) is independent of \(\tau\), then
	\[
	\|\chi(-\tau A)A^{-\mu}\|\le C\tau^\mu,
	\qquad0\le\mu\le1.
	\]
\end{thm}

\begin{proof}
	Fix \(0\le\mu\le1\). Since \(A^{-\mu}\) is defined by the functional calculus,
	\[
	A^{-\mu}=(z\mapsto z^{-\mu})(A).
	\]
	Moreover, by the composition and multiplicativity properties of the
	functional calculus,
	\[
	\chi(-\tau A)A^{-\mu}
	=
	\Bigl(z\mapsto z^{-\mu}\chi(-\tau z)\Bigr)(A).
	\]
	Define
	\[
	f_{\mu,\tau}(z):=z^{-\mu}\chi(-\tau z),
	\qquad z\in\Omega_\theta.
	\]
	Since \(\chi\) is holomorphic on a sector containing \(-\tau\Omega_\theta\),
	the function \(f_{\mu,\tau}\) is holomorphic on \(\Omega_\theta\).
	By Assumption~\ref{ass1},
	\[   
	\|\chi(-\tau A)A^{-\mu}\|
	=
	\|f_{\mu,\tau}(A)\|
	\le
	C_A\|f_{\mu,\tau}\|_{H^\infty(\Omega_\theta)}.
	\]
	Using the change of variables \(w=-\tau z\), we obtain
	\[
	\begin{aligned}
		\|f_{\mu,\tau}\|_{H^\infty(\Omega_\theta)}
		&=
		\sup_{z\in\Omega_\theta}
		|z^{-\mu}\chi(-\tau z)| \\
		&=
		\sup_{w\in-\tau\Omega_\theta}
		\left|\left(-\frac{w}{\tau}\right)^{-\mu}\chi(w)\right| \\
		&=
		\tau^\mu
		\sup_{w\in-\tau\Omega_\theta}
		|w^{-\mu}\chi(w)| \\
		&\le
		C_\chi\tau^\mu.
	\end{aligned}
	\]
	Therefore,
	\[
	\|\chi(-\tau A)A^{-\mu}\|
	\le
	C_A C_\chi\tau^\mu.
	\]
	Setting \(C:=C_A C_\chi\) yields the asserted estimate.
\end{proof}

In addition, the following smoothing bounds are used repeatedly \cite{HOCHBRUCK2005323}.

\begin{lem}\label{lem1}
	Under Assumption~\ref{ass1}, for every $\gamma\ge 0$ there is a constant $C=C(\gamma,T)$ such that
	\begin{equation}\label{parabolicsmoothing1}
		\|\mathrm{e}^{-tA}\| + \|t^{\gamma}A^{\gamma}\mathrm{e}^{-tA}\| \le C, \qquad t\in[0,T].
	\end{equation}
\end{lem}

\begin{assumption}\label{ass2}
	Let $0<\gamma<1$ and let $B:\mathcal{D}(B)\to X$ be a closed linear operator such that $\mathcal{D}(A^{\gamma})\subset\mathcal{D}(B)$ and
	\begin{equation}\label{eq99}
		\|Bx\|\le \varepsilon\|A^{\gamma}x\| + K(\varepsilon)\|x\| \qquad \text{for all } x\in\mathcal{D}(A^{\gamma}),
	\end{equation}
	for sufficiently small $\varepsilon>0$, where $K(\varepsilon)\ge 0$. Assumption~\ref{ass2} implies \cite{pazy1983semigroups} 

\begin{equation}\label{eq14}
	\|B A^{-\gamma}\|\le C.
\end{equation}

	In addition, we assume that
	\begin{equation}\label{eq15}
		\|A^{-\gamma}B\|\le C.
	\end{equation}
\end{assumption}

\begin{lem}[\cite{Hau2025}]\label{lem3}
	Under Assumptions~\ref{ass1} and \ref{ass2}, let $\phi(-\tau A)$ be an arbitrary linear combination of the functions $\varphi_k(-\tau A)$. Then
	\begin{subequations}
		\begin{align}
			\|\phi(-\tau A)B\| &\le C\tau^{-\gamma}, \label{eq68a}\\
			\|\mathrm{e}^{-tA}\phi(-\tau A)B\| &\le Ct^{-\gamma}. \label{eq68b}
		\end{align}
	\end{subequations}
\end{lem}

Under Assumptions~\ref{ass1} and \ref{ass2}, the operator $A-B$ is again sectorial \cite{pazy1983semigroups}. Hence Lemma~\ref{lem1} applies to $\mathrm{e}^{-t(A-B)}$ and yields bounds on the derivatives of $u(t)=\mathrm{e}^{-t(A-B)}u_0$. Throughout this work, we assume that $u$ possesses sufficient regularity and that both $u$ and its derivatives belong to the domain of $B$. Concrete examples of operators satisfying Assumptions~\ref{ass1} and \ref{ass2} will be given later.

\section{Exponential Runge--Kutta methods}
\label{recallEXP}

Let $\tau>0$ denote the time step size, and set
$t_n=n\tau$ and $t_{n+1}=t_n+\tau$.
We write $u_n$ for the numerical approximation to
$u(t_n)$ and $U_{ni}$ for the internal approximation
associated with the intermediate time $t_n+c_i\tau$.
An explicit exponential Runge--Kutta method with $s$ stages
applied to \eqref{eq:model} is defined by
\begin{subequations}\label{eq20}
	\begin{equation}\label{eq20a}
		U_{ni}
		=
		\mathrm{e}^{-c_i\tau A}u_n
		+
		\tau\sum_{j=1}^{i-1}
		a_{ij}(-\tau A)BU_{nj},
		\qquad 1\leq i\leq s,
	\end{equation}
	\begin{equation}\label{eq20b}
		u_{n+1}
		=
		\mathrm{e}^{-\tau A}u_n
		+
		\tau\sum_{i=1}^{s}
		b_i(-\tau A)BU_{ni}.
	\end{equation}
\end{subequations}
We take $c_1=0$ and $U_{n1}=u_n$, while
$c_i\in[0,1]$ for $2\leq i\leq s$.
The operator-valued coefficient functions
$a_{ij}(-\tau A)$ and $b_i(-\tau A)$ are constructed
from linear combinations of $\varphi$-functions,
where
\[
\varphi_k(z)
=
\int_0^1
e^{(1-\theta)z}
\frac{\theta^{k-1}}{(k-1)!}\,d\theta,
\qquad k\geq1.
\]
Their values at $-\tau A$ are well defined by
Lemma \ref{parabolicsmoothing1}.
Furthermore, these coefficient functions inherit the
smoothing property of the analytic semigroup:
\begin{equation}\label{eq19}
	\|\phi(-tA)\|
	+
	\|t^\eta A^\eta\phi(-tA)\|
	\leq C,
	\qquad 0\leq\eta\leq1,
\end{equation}
where $\phi(-tA)$ represents any coefficient function
appearing in the numerical scheme.

We next specify the class of fourth-order methods
considered in this work. For this purpose, introduce,
for $2\leq i\leq s$, the stage functions
\begin{equation}\label{eq:psi2_fourth}
	\psi_{2,i}(-\tau A)
	=
	c_i^2\varphi_2(-c_i\tau A)
	-
	\sum_{j=2}^{i-1}
	c_j a_{ij}(-\tau A),
\end{equation}
and
\begin{equation}\label{eq:psi3_fourth}
	\psi_{3,i}(-\tau A)
	=
	c_i^3\varphi_3(-c_i\tau A)
	-
	\frac{1}{2}
	\sum_{j=2}^{i-1}
	c_j^2a_{ij}(-\tau A).
\end{equation}
In terms of these functions, the coefficient functions
of the method are required to satisfy the following
conditions:
{\allowdisplaybreaks
\begin{subequations}\label{eq:fourth_order_conditions}
	\begin{align}
		&\sum_{i=1}^{s}b_i(-\tau A)
		=
		\varphi_1(-\tau A),
		\label{eq:fourth_order_conditions_a}\\
		&\sum_{j=1}^{i-1}a_{ij}(-\tau A)
		=
		c_i\varphi_1(-c_i\tau A),
		\qquad 2\leq i\leq s,
		\label{eq:fourth_order_conditions_b}\\
		&\sum_{i=2}^{s}c_i b_i(-\tau A)
		=
		\varphi_2(-\tau A),
		\label{eq:fourth_order_conditions_c}\\
		&\sum_{i=2}^{s}c_i^2 b_i(-\tau A)
		=
		2\varphi_3(-\tau A),
		\label{eq:fourth_order_conditions_d}\\
		&\sum_{i=2}^{s}
		b_i(-\tau A)J\psi_{2,i}(-\tau A)
		=
		0,
		\label{eq:fourth_order_conditions_e}\\
		&\sum_{i=2}^{s}
		\frac{c_i^3}{3!}b_i(0)
		=
		\varphi_4(0)
		=
		\frac{1}{4!},
		\label{eq:fourth_order_conditions_f}\\
		&\sum_{i=2}^{s}
		b_i(-\tau A)J\psi_{3,i}(-\tau A)
		=
		0,
		\label{eq:fourth_order_conditions_g}\\
		&\sum_{i=3}^{s}
		b_i(-\tau A)J
		\left(
		\sum_{j=2}^{i-1}
		a_{ij}(-\tau A)J\psi_{2,j}(-\tau A)
		\right)
		=
		0,
		\label{eq:fourth_order_conditions_h}\\
		&\sum_{i=2}^{s}
		c_i b_i(-\tau A)K\psi_{2,i}(-\tau A)
		=
		0.
		\label{eq:fourth_order_conditions_i}
	\end{align}
\end{subequations}
}
Here, $J$ and $K$ are arbitrary operators of compatible
dimensions.

The first two relations in
\eqref{eq:fourth_order_conditions} ensure the consistency
of the stages and the final update. The remaining
identities encode the additional requirements associated
with fourth-order accuracy in the stiff setting.
They correspond to the fourth-order stiff order
conditions introduced in \cite{Luan2021}, expressed in
the present exponential Runge--Kutta formulation.
An important feature of the class is that the condition
\eqref{eq:fourth_order_conditions_f} is required only
at the origin, whereas all the other stiff order
conditions are imposed in their full operator-dependent
form. This particular relaxation is sufficient for
the fourth-order convergence result established in
the semilinear parabolic setting in \cite{Luan2021}.

The class described by
\eqref{eq:fourth_order_conditions} is not restricted
to a particular choice of the number of stages or
the values of the intermediate nodes.
For example, the six-stage method $\mathrm{expRK4s6}$
constructed in \cite{Luan2021} belongs to this class.
In the following, however, we consider general
$s$-stage methods whose coefficients satisfy the
above relations, rather than restricting the analysis
to this particular construction.

The convergence properties of these methods have
previously been investigated for semilinear parabolic
equations with sufficiently regular nonlinearities.
The situation considered here is different: the
linear model problem \eqref{eq:model} involves the
operator $B$ directly, and $B$ is not assumed to be
bounded or to commute with $A$.
Thus, the usual convergence arguments for bounded
nonlinear perturbations do not immediately provide
the required estimates in our setting.

The purpose of this paper is to extend the convergence
analysis to this class of problems. Under suitable
assumptions on the operators and the regularity of
the initial data, we establish fourth-order convergence
for explicit exponential Runge--Kutta methods satisfying
\eqref{eq:fourth_order_conditions}.
A central part of the analysis is to examine the
interplay between the stiff order conditions and the
non-commutativity of $A$ and $B$.
In particular, we identify the cancellations that
remove the defect terms responsible for order reduction
and establish bounds for the remaining contributions.
This provides a theoretical explanation for the
preservation of fourth-order accuracy in the presence
of unbounded, non-commuting operators.

\section{Convergence analysis}
\label{erroranalysis}

This section is devoted to the convergence analysis of the
six-stage exponential Runge--Kutta methods introduced above,
applied to the linear problem \eqref{eq:model}.
Throughout, $C=C(T,\widetilde{\Gamma},\zeta,\gamma)$
denotes a generic positive constant, which may change from
one occurrence to another. We point out that
$\widetilde{\Gamma}$ and $\gamma$ arise from different
operator estimates and are independent parameters.
Moreover, $\zeta>0$ is fixed and sufficiently small.

\subsection{Error recursion and defect representation}
\label{errorrecursion}

The starting point of our analysis is the
variation-of-constants formula, which gives, for
$0\leq\theta\leq1$,
\begin{equation}\label{eq91}
	u(t_n+\theta\tau)
	=
	\mathrm{e}^{-\theta\tau A}u(t_n)
	+
	\int_0^{\theta\tau}
	\mathrm{e}^{-(\theta\tau-\xi)A}
	Bu(t_n+\xi)\,\mathrm{d}\xi.
\end{equation}
To derive the local defect representations, we employ
Taylor's formula in the form
\begin{equation}\label{eq90}
	u(t_n+\xi)
	=
	\sum_{j=1}^{4}
	\frac{\xi^{j-1}}{(j-1)!}
	u^{(j-1)}(t_n)
	+
	\int_0^\xi
	\frac{(\xi-\sigma)^3}{3!}
	u^{(4)}(t_n+\sigma)\,\mathrm{d}\sigma,
	\qquad 0\leq\xi\leq\tau.
\end{equation}
Combining \eqref{eq90} with \eqref{eq91} and evaluating
the integrals associated with the Taylor polynomial yields
\begin{equation}\label{eq87}
	\begin{aligned}
		u(t_n+c_i\tau)
		={}&
		\mathrm{e}^{-c_i\tau A}u(t_n)
		+
		\sum_{j=1}^{4}
		(c_i\tau)^j
		\varphi_j(-c_i\tau A)Bu^{(j-1)}(t_n)
		\\
		&+
		\int_0^{c_i\tau}
		\mathrm{e}^{-(c_i\tau-\xi)A}
		\int_0^\xi
		\frac{(\xi-\sigma)^3}{3!}
		Bu^{(4)}(t_n+\sigma)
		\,\mathrm{d}\sigma\,\mathrm{d}\xi.
	\end{aligned}
\end{equation}

We now apply the numerical scheme to the exact solution.
For the internal stages and the final update, respectively,
we obtain
\begin{subequations}\label{eq84}
	\begin{equation}\label{eq84a}
		u(t_n+c_i\tau)
		=
		\mathrm{e}^{-c_i\tau A}u(t_n)
		+
		\tau\sum_{j=1}^{i-1}
		a_{ij}(-\tau A)Bu(t_n+c_j\tau)
		+
		\Delta_{ni},
		\qquad 2\leq i\leq6,
	\end{equation}
	\begin{equation}\label{eq84b}
		u(t_{n+1})
		=
		\mathrm{e}^{-\tau A}u(t_n)
		+
		\tau\sum_{i=1}^{6}
		b_i(-\tau A)Bu(t_n+c_i\tau)
		+
		\delta_{n+1}.
	\end{equation}
\end{subequations}
Here, $\Delta_{ni}$ and $\delta_{n+1}$ represent the
defects in the internal stages and in the final update,
respectively.

For a more explicit representation of the stage defects,
we substitute \eqref{eq90} into the terms
$u(t_n+c_j\tau)$ appearing in \eqref{eq84a}. This gives
\begin{equation}\label{eq86}
	\begin{aligned}
		u(t_n+c_i\tau)
		={}&
		\mathrm{e}^{-c_i\tau A}u(t_n)
		\\
		&+
		\tau\sum_{k=1}^{i-1}
		a_{ik}(-\tau A)
		\sum_{j=1}^{4}
		\frac{(c_k\tau)^{j-1}}{(j-1)!}
		Bu^{(j-1)}(t_n)
		\\
		&+
		\tau\sum_{k=1}^{i-1}
		a_{ik}(-\tau A)
		\int_0^{c_k\tau}
		\frac{(c_k\tau-\sigma)^3}{3!}
		Bu^{(4)}(t_n+\sigma)\,\mathrm{d}\sigma
		+
		\Delta_{ni}.
	\end{aligned}
\end{equation}
Similarly, inserting the Taylor expansion into
\eqref{eq84b} gives
\begin{equation}\label{eq:fourth_update_expansion}
	\begin{aligned}
		u(t_{n+1})
		={}&
		\mathrm{e}^{-\tau A}u(t_n)
		+
		\tau\sum_{i=1}^{6}
		b_i(-\tau A)
		\sum_{j=1}^{4}
		\frac{(c_i\tau)^{j-1}}{(j-1)!}
		Bu^{(j-1)}(t_n)
		\\
		&+
		\tau\sum_{i=1}^{6}
		b_i(-\tau A)
		\int_0^{c_i\tau}
		\frac{(c_i\tau-\sigma)^3}{3!}
		Bu^{(4)}(t_n+\sigma)\,\mathrm{d}\sigma
		+
		\delta_{n+1}.
	\end{aligned}
\end{equation}

Consequently, the defects $\Delta_{ni}$ and
$\delta_{n+1}$ are obtained by comparing
\eqref{eq86} with \eqref{eq87}, and
\eqref{eq:fourth_update_expansion} with the corresponding
exact expansion at $t_{n+1}$.
Their representations will be used below to identify
the contributions that vanish by virtue of the
fourth-order stiff order conditions and to estimate
the remaining terms.

Finally, define the errors at the time nodes and at
the internal stages by
\[
e_n=u_n-u(t_n),
\qquad
E_{ni}=U_{ni}-u(t_n+c_i\tau).
\]
Subtracting the exact relations \eqref{eq84} from
the numerical scheme yields
\begin{subequations}\label{eq85}
	\begin{equation}\label{eq85a}
		E_{ni}
		=
		\mathrm{e}^{-c_i\tau A}e_n
		+
		\tau\sum_{j=1}^{i-1}
		a_{ij}(-\tau A)BE_{nj}
		-
		\Delta_{ni},
		\qquad 2\leq i\leq6,
	\end{equation}
	\begin{equation}\label{eq85b}
		e_{n+1}
		=
		\mathrm{e}^{-\tau A}e_n
		+
		\tau\sum_{i=1}^{6}
		b_i(-\tau A)BE_{ni}
		-
		\delta_{n+1}.
	\end{equation}
\end{subequations}
These relations form the basis of the subsequent
convergence analysis. In particular, the fourth-order
stiff order conditions will be used to control the
contributions of the stage defects and the final-update
defect, taking into account the non-commutativity
and unboundedness of the operators $A$ and $B$.

\subsection{Convergence results for fourth-order methods}

We are now ready to investigate the convergence of the
six-stage fourth-order exponential Runge--Kutta methods
satisfying the order conditions
\eqref{eq:fourth_order_conditions}. Throughout this section,
we assume that
\[
u_0\in\mathcal{D}((A-B)^3).
\]

\begin{rem}
	As in the third-order analysis, we distinguish between
	the initial time interval and subsequent time steps.
	Indeed, the regularity available at $t=0$ allows us to
	use the Taylor expansion
	\begin{equation}\label{eq:fourth_initial_taylor}
		u(\xi)
		=
		u(0)+\xi u'(0)+\frac{\xi^2}{2}u''(0)
		+
		\int_0^\xi
		\frac{(\xi-\sigma)^2}{2}
		u'''(\sigma)\,\mathrm{d}\sigma,
	\end{equation}
	whereas, for $t_n>0$, we employ
	\begin{equation}\label{eq:fourth_positive_taylor}
		\begin{aligned}
			u(t_n+\xi)
			={}&u(t_n)+\xi u'(t_n)
			+\frac{\xi^2}{2}u''(t_n)
			+\frac{\xi^3}{3!}u'''(t_n)
			\\
			&+
			\int_0^\xi
			\frac{(\xi-\sigma)^3}{3!}
			u^{(4)}(t_n+\sigma)\,\mathrm{d}\sigma.
		\end{aligned}
	\end{equation}
	Recall that the solution of \eqref{eq:model} satisfies
	\[
	u(t)=\mathrm{e}^{-t(A-B)}u_0.
	\]
	In particular, its third- and fourth-order derivatives
	are given by
	\begin{equation}\label{eq:fourth_derivatives}
		\begin{aligned}
			u'''(t)
			&=-\mathrm{e}^{-t(A-B)}(A-B)^3u_0,\\
			u^{(4)}(t)
			&=\mathrm{e}^{-t(A-B)}(A-B)^4u_0.
		\end{aligned}
	\end{equation}
\end{rem}

The following lemma provides the regularity estimates
needed in the error analysis.

\begin{lem}\label{lem:fourth_regularity}
	Under Assumption \ref{ass1}, for $t>0$ and
	$u_0\in\mathcal{D}((A-B)^3)$, we have
	\begin{subequations}
		\begin{equation}\label{eq:fourth_derivative_bound}
			\|u'''(t)\|\leq C,
		\end{equation}
		\begin{equation}\label{eq:fourth_derivative_smoothing}
			\|u^{(4)}(t)\|\leq Ct^{-1}.
		\end{equation}
	\end{subequations}
\end{lem}

\begin{proof}
	The first estimate follows from the boundedness of the
	analytic semigroup generated by $-(A-B)$ and the
	assumption $u_0\in\mathcal{D}((A-B)^3)$.
	For the second estimate, we write
	\[
	u^{(4)}(t)
	=
	(A-B)\mathrm{e}^{-t(A-B)}(A-B)^3u_0
	\]
	and apply the standard analytic semigroup smoothing
	estimate.
\end{proof}

Following the error recursion derived in
Section \ref{errorrecursion}, we substitute the defect
representations into \eqref{eq85} and successively
eliminate the internal stage errors. Using the fourth-order
stiff order conditions \eqref{eq:fourth_order_conditions},
we obtain
\begin{equation}\label{eq:fourth_global_recursion}
	\begin{aligned}
		e_{n+1}
		={}&
		\mathrm{e}^{-\tau A}e_n
		+\tau\mathcal{N}e_n
		+\tau^2\mathcal{M}e_n
		+\tau^3\mathcal{K}e_n
		+\tau^4\mathcal{P}e_n
		\\
		&+\tau^5\mathcal{Q}e_n
		+\tau^6\mathcal{S}e_n
		+\mathcal{T}_{n+1}
		+\mathcal{R}_{n+1}
		+\widehat{\mathcal{R}}_{n+1},
	\end{aligned}
\end{equation}
where the operators independent of $n$ are given by
\begin{equation*}
	\begin{aligned}
		\mathcal{N}
		={}&
		\sum_{i=1}^{6}
		b_i(-\tau A)B\mathrm{e}^{-c_i\tau A},
		\\[1mm]
		\mathcal{M}
		={}&
		\sum_{i=2}^{6}\sum_{j=1}^{i-1}
		b_i(-\tau A)B
		a_{ij}(-\tau A)B
		\mathrm{e}^{-c_j\tau A},
		\\[1mm]
		\mathcal{K}
		={}&
		\sum_{i=3}^{6}\sum_{j=2}^{i-1}\sum_{k=1}^{j-1}
		b_i(-\tau A)B
		a_{ij}(-\tau A)B
		a_{jk}(-\tau A)B
		\mathrm{e}^{-c_k\tau A},
		\\[1mm]
		\mathcal{P}
		={}&
		\sum_{i=4}^{6}\sum_{j=3}^{i-1}
		\sum_{k=2}^{j-1}\sum_{\ell=1}^{k-1}
		b_i(-\tau A)B
		a_{ij}(-\tau A)B
		a_{jk}(-\tau A)B
		a_{k\ell}(-\tau A)B
		\mathrm{e}^{-c_\ell\tau A},
		\\[1mm]
		\mathcal{Q}
		={}&
		\sum_{i=5}^{6}\sum_{j=4}^{i-1}
		\sum_{k=3}^{j-1}\sum_{\ell=2}^{k-1}
		\sum_{r=1}^{\ell-1}
		b_i(-\tau A)B
		a_{ij}(-\tau A)B
		a_{jk}(-\tau A)B
		a_{k\ell}(-\tau A)B
		a_{\ell r}(-\tau A)B
		\mathrm{e}^{-c_r\tau A},
		\\[1mm]
		\mathcal{S}
		={}&
		b_6(-\tau A)B
		a_{65}(-\tau A)B
		a_{54}(-\tau A)B
		a_{43}(-\tau A)B
		a_{32}(-\tau A)B
		a_{21}(-\tau A)B.
	\end{aligned}
\end{equation*}

The contribution of the weakened fourth-order quadrature
condition is contained in
\begin{equation}\label{eq:fourth_T}
	\mathcal{T}_{n+1}
	=
	\tau^4
	\left(
	\sum_{i=2}^{6}\frac{c_i^3}{3!}b_i(-\tau A)
	-
	\varphi_4(-\tau A)
	\right)
	Bu'''(t_n),
	\qquad n\geq1.
\end{equation}
Indeed, the fourth-order quadrature condition
\eqref{eq:fourth_order_conditions_f} implies that
\[
\sum_{i=2}^{6}\frac{c_i^3}{3!}b_i(0)
-
\varphi_4(0)=0.
\]
Thus, the operator-valued coefficient in
\eqref{eq:fourth_T} vanishes at the origin, a property
that will be exploited in estimating this contribution.

The terms arising from the internal stage defects are
collected in
\begin{equation}\label{eq:remainder_after_cancellation}
	\begin{aligned}
		\mathcal R_{n+1}
		={}&
		-\tau\sum_{i=2}^{6}
		b_i(-\tau A)B\Theta_{ni}
		\\
		&-\tau^2\sum_{i=3}^{6}\sum_{j=2}^{i-1}
		b_i(-\tau A)B a_{ij}(-\tau A)B
		\left(
		\tau^3\psi_{3,j}(-\tau A)Bu''(t_n)
		+\Theta_{nj}
		\right)
		\\
		&-\tau^3\sum_{i=4}^{6}\sum_{j=3}^{i-1}
		\sum_{k=2}^{j-1}
		b_i(-\tau A)B a_{ij}(-\tau A)B
		a_{jk}(-\tau A)B
		\\
		&\qquad\qquad\times
		\left(
		\tau^2\psi_{2,k}(-\tau A)Bu'(t_n)
		+\tau^3\psi_{3,k}(-\tau A)Bu''(t_n)
		+\Theta_{nk}
		\right)
		\\
		&-\tau^4\sum_{i=5}^{6}\sum_{j=4}^{i-1}
		\sum_{k=3}^{j-1}\sum_{\ell=2}^{k-1}
		b_i(-\tau A)B a_{ij}(-\tau A)B
		a_{jk}(-\tau A)B a_{k\ell}(-\tau A)B
		\\
		&\qquad\qquad\times
		\left(
		\tau^2\psi_{2,\ell}(-\tau A)Bu'(t_n)
		+\tau^3\psi_{3,\ell}(-\tau A)Bu''(t_n)
		+\Theta_{n\ell}
		\right)
		\\
		&-\tau^5 b_6(-\tau A)B a_{65}(-\tau A)B
		a_{54}(-\tau A)B a_{43}(-\tau A)B
		a_{32}(-\tau A)B
		\\
		&\qquad\qquad\times
		\left(
		\tau^2\psi_{2,2}(-\tau A)Bu'(t_n)
		+\tau^3\psi_{3,2}(-\tau A)Bu''(t_n)
		+\Theta_{n2}
		\right),
	\end{aligned}
\end{equation}
where \begin{equation}
	\begin{aligned}
		\Theta_{ni}
		={}&
		\int_0^{c_i\tau}
		\mathrm{e}^{-(c_i\tau-\sigma)A}B
		\int_0^\sigma
		\frac{(\sigma-\rho)^2}{2}
		u'''(t_n+\rho)\,\mathrm{d}\rho
		\,\mathrm{d}\sigma
		\\
		&-\tau\sum_{j=2}^{i-1}
		a_{ij}(-\tau A)B
		\int_0^{c_j\tau}
		\frac{(c_j\tau-\rho)^2}{2}
		u'''(t_n+\rho)\,\mathrm{d}\rho.
	\end{aligned}
\end{equation}

For $n\geq1$, $\widehat{\mathcal{R}}_{n+1}$ is given by
\begin{equation}\label{eq:fourth_Rhat}
	\begin{aligned}
		\widehat{\mathcal{R}}_{n+1}
		={}&
		-\int_0^\tau
		\mathrm{e}^{-(\tau-\xi)A}
		\int_0^\xi
		\frac{(\xi-\sigma)^3}{3!}
		Bu^{(4)}(t_n+\sigma)
		\,\mathrm{d}\sigma\,\mathrm{d}\xi
		\\
		&+
		\tau\sum_{i=2}^{6}
		b_i(-\tau A)B
		\int_0^{c_i\tau}
		\frac{(c_i\tau-\sigma)^3}{3!}
		u^{(4)}(t_n+\sigma)\,\mathrm{d}\sigma.
	\end{aligned}
\end{equation}
For the first time step, the expansion
\eqref{eq:fourth_initial_taylor} is used instead, and
the corresponding remainder is
\begin{equation}\label{eq:fourth_Rhat_initial}
	\begin{aligned}
		\widehat{\mathcal{R}}_1
		={}&
		-\int_0^\tau
		\mathrm{e}^{-(\tau-\xi)A}
		\int_0^\xi
		\frac{(\xi-\sigma)^2}{2}
		Bu'''(\sigma)\,\mathrm{d}\sigma\,\mathrm{d}\xi
		\\
		&+
		\tau\sum_{i=2}^{6}
		b_i(-\tau A)B
		\int_0^{c_i\tau}
		\frac{(c_i\tau-\sigma)^2}{2}
		u'''(\sigma)\,\mathrm{d}\sigma.
	\end{aligned}
\end{equation}

Solving \eqref{eq:fourth_global_recursion} gives
\begin{equation}\label{eq:fourth_error_solution}
	\begin{aligned}
		e_n
		={}&
		\tau\sum_{j=0}^{n-1}
		\mathrm{e}^{-(n-j-1)\tau A}\mathcal{N}e_j
		+
		\tau^2\sum_{j=0}^{n-1}
		\mathrm{e}^{-(n-j-1)\tau A}\mathcal{M}e_j
		\\
		&+
		\tau^3\sum_{j=0}^{n-1}
		\mathrm{e}^{-(n-j-1)\tau A}\mathcal{K}e_j
		+
		\tau^4\sum_{j=0}^{n-1}
		\mathrm{e}^{-(n-j-1)\tau A}\mathcal{P}e_j
		\\
		&+
		\tau^5\sum_{j=0}^{n-1}
		\mathrm{e}^{-(n-j-1)\tau A}\mathcal{Q}e_j
		+
		\tau^6\sum_{j=0}^{n-1}
		\mathrm{e}^{-(n-j-1)\tau A}\mathcal{S}e_j
		\\
		&+
		\sum_{j=0}^{n-1}
		\mathrm{e}^{-(n-j-1)\tau A}
		\left(
		\mathcal{T}_{j+1}
		+\mathcal{R}_{j+1}
		+\widehat{\mathcal{R}}_{j+1}
		\right).
	\end{aligned}
\end{equation}
Consequently, the convergence estimate reduces to
bounding the six error propagation contributions and
the three defect contributions in
\eqref{eq:fourth_error_solution}.
The analytic semigroup estimates, the regularity bounds
of Lemma \ref{lem:fourth_regularity}, and the fourth-order
stiff order conditions will be used to obtain these
bounds.

\subsection{Preliminary error bounds}
\label{error-analysis}

We first estimate the terms arising from the propagation of the
previous errors. By the definition of $\mathcal{N}$ and
\eqref{eq68a}--\eqref{eq68b}, we have
\[
\left\|\mathrm{e}^{-sA}\mathcal{N}\right\|
\leq Cs^{-\gamma},
\qquad
\|\mathcal{N}\|\leq C\tau^{-\gamma}.
\]
Consequently, separating the last term in the sum, we obtain
\begin{align}
	&\left\|
	\tau\sum_{j=0}^{n-1}
	\mathrm{e}^{-(n-j-1)\tau A}\mathcal{N}e_j
	\right\|
	\nonumber\\
	&\leq
	C\tau\sum_{j=0}^{n-2}
	\left\|
	\mathrm{e}^{-(n-j-1)\tau A}\mathcal{N}
	\right\|\|e_j\|
	+
	C\tau\|\mathcal{N}\|\|e_{n-1}\|
	\nonumber\\
	&\leq
	C\tau\sum_{j=0}^{n-2}
	t_{n-j-1}^{-\gamma}\|e_j\|
	+
	C\tau^{1-\gamma}\|e_{n-1}\|.
	\label{eq:fourth_N_bound}
\end{align}

The estimate for $\mathcal{M}$ follows in the same way. Indeed, each
term in $\mathcal{M}$ contains, in addition to the corresponding
factors in $\mathcal{N}$, one extra stage factor of the form
$\tau a_{ik}(-\tau A)B$. By \eqref{eq68a}, this factor satisfies
\[
\tau\|a_{ik}(-\tau A)B\|
\leq C\tau^{1-\gamma}.
\]
Thus, we obtain
\begin{align}
	&\left\|
	\tau^2\sum_{j=0}^{n-1}
	\mathrm{e}^{-(n-j-1)\tau A}\mathcal{M}e_j
	\right\|
	\nonumber\\
	&\leq
	C\tau^{1-\gamma}
	\left(
	\tau\sum_{j=0}^{n-2}
	t_{n-j-1}^{-\gamma}\|e_j\|
	+
	\tau^{1-\gamma}\|e_{n-1}\|
	\right)
	\nonumber\\
	&\leq
	C\tau^{2-\gamma}
	\sum_{j=0}^{n-2}
	t_{n-j-1}^{-\gamma}\|e_j\|
	+
	C\tau^{2-2\gamma}\|e_{n-1}\|.
	\label{eq:fourth_M_bound}
\end{align}

The same argument applies to $\mathcal{K}$, $\mathcal{P}$,
$\mathcal{Q}$, and $\mathcal{S}$. The only difference is that their
terms contain, respectively, two, three, four, and five additional
stage factors of the form $\tau a_{ik}(-\tau A)B$. Hence, each
successive contribution gains another factor $\tau^{1-\gamma}$.
In particular, for $r=2,3,4,5$, corresponding respectively to
$\mathcal{K}$, $\mathcal{P}$, $\mathcal{Q}$, and $\mathcal{S}$,
we obtain
\begin{align}
	&\left\|
	\tau^{r+1}\sum_{j=0}^{n-1}
	\mathrm{e}^{-(n-j-1)\tau A}
	\mathcal{X}_{r+1}e_j
	\right\|
	\nonumber\\
	&\leq
	C\tau^{r(1-\gamma)}
	\left(
	\tau\sum_{j=0}^{n-2}
	t_{n-j-1}^{-\gamma}\|e_j\|
	+
	\tau^{1-\gamma}\|e_{n-1}\|
	\right),
	\label{eq:fourth_remaining_propagation_bounds}
\end{align}
where $\mathcal{X}_{r+1}$ denotes the corresponding operator
$\mathcal{K}$, $\mathcal{P}$, $\mathcal{Q}$, or $\mathcal{S}$.

We next establish two auxiliary estimates for the terms involving
$B\phi(-\tau A)B$. Recall that $\phi(-\tau A)$ denotes a linear
combination of the functions $\varphi_k(-\tau A)$.

\begin{lem}
	\label{lem:fourth_BphiB_uprime}
	Let $\gamma\geq\widetilde{\Gamma}$, $t\in(0,T]$, and assume that
	\[
	\|A^{\widetilde{\Gamma}}BA^{-2}\|\leq C,
	\qquad
	u_0\in\mathcal{D}((A-B)^3).
	\]
	Then
	\begin{align}
		\|B\phi(-\tau A)Bu'(t)\|
		&\leq C\tau^{\widetilde{\Gamma}-\gamma},
		\label{eq:fourth_BphiB_uprime}
		\\
		\|B\phi(-\tau A)Bu'''(t)\|
		&\leq Ct^{-2}\tau^{\widetilde{\Gamma}-\gamma}.
		\label{eq:fourth_BphiB_uprimeprimeprime}
	\end{align}
\end{lem}
\begin{proof}
	By \eqref{eq14}, \eqref{eq19}, and the analytic smoothing
	property, we obtain
	\begin{align*}
		\|B\phi(-\tau A)Bu'(t)\|
		&\leq
		\|BA^{-\gamma}\|
		\|\phi(-\tau A)A^{\gamma-\widetilde{\Gamma}}\|
		\|A^{\widetilde{\Gamma}}BA^{-2}\|
		\\
		&\qquad\times
		\|A^2(A-B)^{-2}\|
		\|u'''(t)\|
		\\
		&\leq
		C\tau^{\widetilde{\Gamma}-\gamma}.
	\end{align*}
	Similarly, applying the same factorization to
	$B\phi(-\tau A)Bu'''(t)$, we obtain
	\begin{align*}
		\|B\phi(-\tau A)Bu'''(t)\|
		&\leq
		\|BA^{-\gamma}\|
		\|\phi(-\tau A)A^{\gamma-\widetilde{\Gamma}}\|
		\|A^{\widetilde{\Gamma}}BA^{-2}\|
		\\
		&\qquad\times
		\|A^2(A-B)^{-2}\|
		\|u^{(5)}(t)\|
		\\
		&\leq
		Ct^{-2}\tau^{\widetilde{\Gamma}-\gamma}.
	\end{align*}
	Here, we used
	\[
	u'''(t)
	=-(A-B)^3\mathrm{e}^{-t(A-B)}u_0,
	\qquad
	\|u'''(t)\|\leq C,
	\]
	and
	\[
	u^{(5)}(t)
	=-(A-B)^2\mathrm{e}^{-t(A-B)}(A-B)^3u_0,
	\qquad
	\|u^{(5)}(t)\|\leq Ct^{-2}.
	\]
	The first estimate follows from the uniform boundedness of the
	semigroup, while the second follows from its analytic smoothing
	property.
\end{proof}

We now estimate the quadrature defect
\[
\mathcal{T}_{n+1}
=
\tau^4
\left(
\sum_{i=2}^{6}\frac{c_i^3}{3!}b_i(-\tau A)
-\varphi_4(-\tau A)
\right)
Bu'''(t_n).
\]
%
For $j=0$, employing Theorem \ref{thm:quadrature-defect}
together with $\|A^{\widetilde{\Gamma}}BA^{-1}\|\leq C$, we obtain
\begin{align*}
	\|\mathcal{T}_{n}\|
	&\leq
	C\tau^4
	\left\|
	\left(
	\sum_{i=2}^{6}\frac{c_i^3}{3!}b_i(-\tau A)
	-\varphi_4(-\tau A)
	\right)A^{-\widetilde{\Gamma}}
	\right\|
	\\
	&\qquad\times
	\|A^{\widetilde{\Gamma}}BA^{-1}\|
	\|A(A-B)^{-1}\|
	\|u^{(4)}(t_{n-1})\|
	\\
	&\leq
	Ct_{n-1}^{-1}\tau^{4+\widetilde{\Gamma}}.
\end{align*}
For $j=n-1$, recalling the definition of $\mathcal{T}_1$ and
using \eqref{eq68b}, we get
\begin{align*}
	&\left\|
	\mathrm{e}^{-(n-1)\tau A}\mathcal{T}_1
	\right\|
	\\
	&\leq
	C\tau^4
	\left(
	\sum_{i=2}^{6}
	\left\|
	\mathrm{e}^{-(n-1)\tau A}b_i(-\tau A)B
	\right\|
	+
	\left\|
	\mathrm{e}^{-(n-1)\tau A}
	\varphi_4(-\tau A)B
	\right\|
	\right)
	\|u'''(0)\|
	\\
	&\leq
	Ct_{n-1}^{-\gamma}\tau^4.
\end{align*}

For the remaining terms, using Theorem \ref{thm:quadrature-defect}, we deduce
\begin{align*}
	&\left\|
	\sum_{j=1}^{n-2}
	\mathrm{e}^{-j\tau A}\mathcal{T}_{n-j}
	\right\|
	\\
	&\leq
	C\tau^5
	\sum_{j=1}^{n-2}
	\left\|\mathrm{e}^{-j\tau A }A \right\|
	\|A^{\widetilde{\Gamma}}BA^{-1}\|
	\|A(A-B)^{-1}\|
	\|u^{(4)}(t_{n-j-1})\|
	\\
	&\leq
	C\tau^{5}
	\sum_{j=1}^{n-2}
	t_j^{-1}t_{n-j-1}^{-1}
	\\
	&\leq
	C\tau^{5-\zeta}
	\sum_{j=1}^{n-2}
	t_j^{-1+\zeta}t_{n-j-1}^{-1+\zeta}
	\\
	&\leq
	Ct_n^{-1+2\zeta}
	\tau^{4-\zeta}.
\end{align*}
Combining the estimates, we obtain
\begin{equation}
	\label{eq:fourth_T_bound}
	\left\|
	\sum_{j=0}^{n-1}
	\mathrm{e}^{-j\tau A}\mathcal{T}_{n-j}
	\right\|
	\leq
	Ct_n^{-1+2\zeta}
\tau^{4-\zeta}.
\end{equation}

We next estimate
\[
\sum_{j=0}^{n-1}\mathrm{e}^{-j\tau A}\mathcal{R}_{n-j}.
\]
To this end, we first collect the contributions from all internal
stages and apply the stiff order conditions
\eqref{eq:fourth_order_conditions_a}--\eqref{eq:fourth_order_conditions_i}.
These conditions ensure the cancellation of the lower-order terms
in the resulting weighted sums. Consequently, the remaining
contributions involve terms of the form
\begin{equation}\label{eq88}
	\tau^5 b_i(-\tau A)B
	\phi(-\tau A)Bu'''(t_{n-1}),
\end{equation}
where \(\phi(-\tau A)\) denotes a linear combination of
\(\varphi\)-functions. The other contributions arise from further
substitutions of the internal stage expansions and contain
additional factors of the form \(\tau a_{ik}(-\tau A)B\).
These factors are uniformly bounded by \eqref{eq68a}, and the
corresponding terms can be estimated in the same way using the
auxiliary lemma, while preserving the order of the operator
products. Since the number of stages is fixed, it therefore
suffices to estimate \eqref{eq88} as a representative contribution.

For $j=0$, using \eqref{eq68a} and
\eqref{eq:fourth_BphiB_uprimeprimeprime}, we obtain
\begin{align*}
	\|\mathcal{R}_{n}\|
	&\leq
	C\tau^5
	\|b_i(-\tau A)B\|
	\|B\phi(-\tau A)Bu'''(t_{n-1})\|
	\\
	&\leq
	C\tau^{5-\gamma}
	\|B\phi(-\tau A)Bu'''(t_{n-1})\|
	\\
	&\leq
	Ct_{n-1}^{-2}
	\tau^{5+\widetilde{\Gamma}-2\gamma}.
\end{align*}

For $j=n-1$, by the smoothing estimate \eqref{eq68b},
we get 
\begin{align*}
	&\left\|
	\mathrm{e}^{-(n-1)\tau A}\mathcal{R}_1^{(i)}
	\right\|
	\\
	&\quad\leq
	C\tau^5
	\left\|
	\mathrm{e}^{-(n-1)\tau A}
	b_i(-\tau A)B
	\right\|
	\left\|B\phi_i(-\tau A)Bu'''(0)\right\|
	\\
	&\quad\leq
	C t_{n-1}^{-\gamma}\tau^{5-2\gamma}.
\end{align*}

For $1\leq j\leq n-2$, we consider a representative contribution
to the residual of the form
\[
\mathcal{R}_{n-j}^{(i)}
=
\tau^5 b_i(-\tau A)B
\phi_i(-\tau A)Bu'''(t_{n-j-1}),
\]
where $\phi_i(-\tau A)$ is a linear combination of
$\varphi_k(-\tau A)$.

\begin{align*}
	& \left\|
	\sum_{j=1}^{n-2}
	\mathrm{e}^{-j\tau A}\mathcal{R}_{n-j}
	\right\|
	\\
	&\quad\leq
	C\tau^5
	\sum_{j=1}^{n-2}
	\left\|
	\mathrm{e}^{-j\tau A}Ab_i(-\tau A) \right\|  \left\|  A^{-1}BA^{\beta}
	\right\| \\
	& \qquad \times 
	\left\|
	\phi_i(-\tau A) A^{-2\beta} \right\| \left\| A^{\beta} B A^{-1} \right\| \left\| A (A-B)^{-1} \right\| \left\| (A-B)u'''(t_{n-j-1})
	\right\|
	\\
	&\quad\leq
	C\tau^{5-2\zeta}
	\sum_{j=1}^{n-2}
	t_j^{-1}t_{n-j-1}^{-1}
	\\
	&\quad\leq
	C\tau^{5-3\zeta}
	\sum_{j=1}^{n-2}
	t_j^{-1+\zeta}t_{n-j-1}^{-1+\zeta} 
	\\
	&\quad\leq
	Ct_n^{-1+2\zeta} \tau^{4-3\zeta}
	.
\end{align*}

Finally, we estimate
\[
\sum_{j=0}^{n-1}
\mathrm{e}^{-j\tau A}\widehat{\mathcal{R}}_{n-j}.
\]
Notably, it is sufficient to estimate the integral contribution
$\widehat{\mathcal{R}}_n^{[1]}$, since the remaining stage
contributions can be bounded in the same way (see \cite{Hau2025}).

We first estimate the contribution of the first integral.
For $n\geq2$, using \eqref{eq68b} and
$\|u^{(4)}(t)\|\leq Ct^{-1}$, we obtain
\begin{align*}
	\|\widehat{\mathcal{R}}_n^{[1]}\|
	&\leq
	C\int_0^\tau
	(\tau-\xi)^{-\gamma}
	\int_0^\xi
	\frac{(\xi-\sigma)^3}{3!}
	\|u^{(4)}(t_{n-1}+\sigma)\|
	\,\mathrm{d}\sigma\,\mathrm{d}\xi
	\\
	&\leq
	Ct_{n-1}^{-1}
	\int_0^\tau
	(\tau-\xi)^{-\gamma}
	\int_0^\xi
	\frac{(\xi-\sigma)^3}{3!}
	\,\mathrm{d}\sigma\,\mathrm{d}\xi
	\\
	&=
	\frac{C}{4!}t_{n-1}^{-1}
	\int_0^\tau
	(\tau-\xi)^{-\gamma}\xi^4
	\,\mathrm{d}\xi
	\\
	&\leq
	Ct_{n-1}^{-1}\tau^{5-\gamma}.
\end{align*}

For $j=n-1$, using the semigroup property and the smoothing estimate
$\|\mathrm{e}^{-tA}B\|\leq Ct^{-\gamma}$, we obtain
\begin{align*}
	&\left\|
	\mathrm{e}^{-(n-1)\tau A}
	\int_0^\tau
	\mathrm{e}^{-(\tau-\xi)A}
	\int_0^\xi
	\frac{(\xi-\sigma)^2}{2}
	Bu'''(\sigma)\,\mathrm{d}\sigma\,\mathrm{d}\xi
	\right\|
	\\
	&\quad\leq
	C\int_0^\tau
	\left\|\mathrm{e}^{-(n\tau-\xi)A}B\right\|
	\int_0^\xi
	\frac{(\xi-\sigma)^2}{2}
	\|u'''(\sigma)\|\,\mathrm{d}\sigma\,\mathrm{d}\xi
	\\
	&\quad\leq
	Ct_{n-1}^{-\gamma}
	\int_0^\tau
	\int_0^\xi
	\frac{(\xi-\sigma)^2}{2}\,\mathrm{d}\sigma\,\mathrm{d}\xi
	\\
	&\quad\leq
	Ct_{n-1}^{-\gamma}\tau^4.
\end{align*}

For $1\leq j\leq n-2$, we get
\begin{align}
	&\left\|
	\sum_{j=1}^{n-2}
	\mathrm{e}^{-j\tau A}
	\widehat{\mathcal{R}}_{n-j}^{[1]}
	\right\|
	\notag\\
	&\quad\leq
	C\sum_{j=1}^{n-2}
	\int_0^\tau\int_0^\xi
	\left\|
	\mathrm{e}^{-j\tau A}
	\mathrm{e}^{-(\tau-\xi)A}B
	\right\|
	t_{n-j-1}^{-1}
	\frac{(\xi-\sigma)^3}{3!}
	\,\mathrm{d}\sigma\,\mathrm{d}\xi
	\notag\\
	&\quad=
	C\sum_{j=1}^{n-2}
	\int_0^\tau\int_0^\xi
	\left\|
	\mathrm{e}^{-j\tau A}A^\zeta
	\mathrm{e}^{-(\tau-\xi)A}A^{-\zeta}B
	\right\|
	t_{n-j-1}^{-1}
	\frac{(\xi-\sigma)^3}{3!}
	\,\mathrm{d}\sigma\,\mathrm{d}\xi
	\notag\\
	&\quad\leq
	C\sum_{j=1}^{n-2}
	t_j^{-\zeta}
	\int_0^\tau\int_0^\xi
	t_{n-j-1}^{-1}
	\frac{(\xi-\sigma)^3}{3!}
	\,\mathrm{d}\sigma\,\mathrm{d}\xi
	\notag\\
	&\quad\leq
	C\tau^5
	\sum_{j=1}^{n-2}
	t_j^{-\zeta}t_{n-j-1}^{-1}.
	\notag\\
	& \quad 	\leq
	C_Tt_n^{-1+2\zeta}\tau^{4-2\zeta},
\end{align}

 We are ready to state the main result.
\subsection{Error bounds} 

\begin{thm}[Convergence of the six-stage method]
	\label{theo4}
	Let Assumptions~\ref{ass1} and~\ref{ass2} hold. Let $\zeta>0$ be an arbitrarily small
	fixed constant. Suppose that
	\[
	\|A^{-1} BA^{\beta}\|\leq C,
	\qquad
	\|A^{-2}BA^{\widetilde{\Gamma}}\|\leq C,
	\]
	for some $\beta,\widetilde{\Gamma}\in\mathbb R$, with
	$\beta \leq\frac12-\zeta$ and 
	$\widetilde{\Gamma}\leq\frac12-\zeta$.
In addition, assume that
\[
u_0\in\mathcal D((A-B)^3).
\]
The six-stage exponential Runge--Kutta method satisfies the
fourth-order stiff order conditions
\eqref{eq:fourth_order_conditions_a}--\eqref{eq:fourth_order_conditions_i}
in strong form, with the sole exception of condition
\eqref{eq:fourth_order_conditions_f}, which is satisfied only
in the weakened form at \(A=0\). Then, for the numerical solution of the initial value
	problem~\eqref{eq:model}, obtained using the six-stage
	fourth-order exponential Runge--Kutta method in
	\eqref{eq20}, we have
	\[
	\|u_n-u(t_n)\|
	\leq
	C t_n^{-1+2\zeta}
	\tau^{4-2\zeta},
	\qquad \text{uniformly for }1\leq n\leq N,
	\]
	where $t_n=n\tau\leq T$. The constant $C$ is independent
	of $n$ and $\tau$, but may depend on $T$ and the constants
	in the assumptions.
\end{thm}

\begin{proof}
	Let $e_n=u_n-u(t_n)$. By combining the estimates for the
	individual terms in the error recursion, we obtain
	\[
	\begin{aligned}
		\|e_n\|
		\leq{}&
		C\tau\sum_{j=1}^{n-2}
		t_{n-j-1}^{-\frac12}\|e_j\|
		+
		C\tau^{\frac12}\|e_{n-1}\| + 
C t_n^{-1+2\zeta}
\tau^{4-2\zeta}.
	\end{aligned}
	\]
	The asserted estimate now follows by applying a discrete
	Gronwall inequality; see
	\cite[Lemma~2.15]{HO2010}.
\end{proof}

Throughout the analysis, the operators $A$ and $B$ are
considered in an abstract setting. To illustrate the
applicability of the assumptions, we consider the
one-dimensional case in which $A$ is a second-order
differential operator with homogeneous Dirichlet boundary
conditions and $B$ is a first-order differential operator.
In this setting, the boundedness of
$\|A^{-2}BA^{\widetilde{\Gamma}}\|$ follows from
Lemma 4.3 of \cite{HauThird}, while the corresponding bound
for $\|A^{-1}  BA^\beta\|$ is established in
Lemma~4.8 of \cite{Hau2025}. These results provide
concrete examples in which the abstract assumptions
used in our convergence analysis are satisfied.

The choice of a six-stage method is motivated by our interest
in achieving high-order convergence while keeping the number
of internal stages relatively reasonable. Our analysis shows that
the method attains convergence rates arbitrarily close to
$4$ in $L^1(\Omega)$, $L^2(\Omega)$, and
$L^\infty(\Omega)$.

It is possible to recover full fourth-order convergence in
all the considered spaces by introducing an additional
internal stage. Such a seven-stage construction, however,
requires additional degrees of freedom to satisfy the
corresponding stiff order conditions and leads to a more
involved analysis of the resulting operator contributions.
The development and analysis of this alternative approach
are beyond the scope of the present work. Instead, we focus
on establishing the convergence properties of the
six-stage method and quantifying the accuracy that can be
achieved with a smaller number of internal stages.

We also emphasize that the rates derived here are based on
the operator bounds established for the considered Banach
spaces. Whether sharper operator estimates or other
functional settings can lead to improved convergence
properties for the six-stage method is an interesting
question for future investigation.

\section{Numerical investigations}
\label{numericalchap2}

In this section, we present numerical experiments to illustrate
the error bounds established in Theorem~\ref{theo4}. The implementation of exponential Runge--Kutta methods requires
the evaluation of linear combinations of the form
\begin{equation}\label{eq:phi_combination}
	\sum_{j=1}^{q}
	\varphi_j(-\tau A)v_j,
\end{equation}
where \(A\) is the matrix arising from the spatial
discretization and \(v_j\) are given vectors. A standard
approach to evaluating such expressions is the augmented
matrix technique; see, for example, \cite{doi:10.1137/100788860}.
This approach transforms the evaluation of a linear combination
of \(\varphi\)-functions into the action of a single matrix
exponential on an appropriately constructed vector. In our
numerical experiments, the augmented matrix exponential is
evaluated using the built-in MATLAB function \texttt{expm}.

For the representative method satisfying the order conditions stated in Theorem \ref{theo4}, we employ the fourth-order, six-stage exponential Runge--Kutta
method expRK4s6 proposed in \cite{Luan2021}.

For completeness, we recall the formulation of expRK4s6 used
in our experiments. Let
\[
c_2=c_3=\frac12,\qquad
c_4=c_6=\frac13,\qquad
c_5=\frac56.
\]
Writing \(F_h(u)=-(A_h-B_h)u\), define
\[
D_{n,i}:=B_h(U_{n,i}-u_n),
\qquad i=2,\ldots,6.
\]
The internal stages are given by
\begin{align}
	U_{n,2}
	={}&u_n+
	c_2\tau\varphi_1(-c_2\tau A_h)F_h(u_n),
	\label{eq:exprk4s6_stage2}
	\\
	U_{n,k}
	={}&u_n+
	c_k\tau\varphi_1(-c_k\tau A_h)F_h(u_n)
	+\frac{c_k^2}{c_2}\tau
	\varphi_2(-c_k\tau A_h)D_{n,2},
	\quad k=3,4,
	\label{eq:exprk4s6_stage34}
	\\
	U_{n,j}
	={}&u_n+
	c_j\tau\varphi_1(-c_j\tau A_h)F_h(u_n)
	\nonumber\\
	&+\frac{c_j^2\tau}{c_3-c_4}
	\varphi_2(-c_j\tau A_h)
	\left(
	-\frac{c_4}{c_3}D_{n,3}
	+\frac{c_3}{c_4}D_{n,4}
	\right)
	\nonumber\\
	&+\frac{2c_j^3\tau}{c_3-c_4}
	\varphi_3(-c_j\tau A_h)
	\left(
	\frac{1}{c_3}D_{n,3}
	-\frac{1}{c_4}D_{n,4}
	\right),
	\quad j=5,6,
	\label{eq:exprk4s6_stage56}
\end{align}
and the numerical solution is updated according to
\begin{align}
	u_{n+1}
	={}&u_n+
	\tau\varphi_1(-\tau A_h)F_h(u_n)
	\nonumber\\
	&+\frac{\tau}{c_5-c_6}
	\varphi_2(-\tau A_h)
	\left(
	-\frac{c_6}{c_5}D_{n,5}
	+\frac{c_5}{c_6}D_{n,6}
	\right)
	\nonumber\\
	&+\frac{2\tau}{c_5-c_6}
	\varphi_3(-\tau A_h)
	\left(
	\frac{1}{c_5}D_{n,5}
	-\frac{1}{c_6}D_{n,6}
	\right).
	\label{eq:exprk4s6_update}
\end{align}

We consider the following linear advection--diffusion problem
on the interval \((0,1)\):
\begin{equation}\label{problem2}
	\begin{aligned}
		\partial_t u(t,x)-\frac12 \partial_{xx}u(t,x)
		&=\partial_xu(t,x),
		&& (t,x)\in(0,1]\times(0,1),\\
		u(t,0)=u(t,1)&=0,
		&& t\in[0,1],\\
		u(0,x)&=64x^3(1-x)^3,
		&& x\in[0,1].
	\end{aligned}
\end{equation}
Equivalently, this problem can be written in the abstract
form
\[
u'(t)+Au(t)=Bu(t),
\]
where
\[
A=-\frac12 \partial_{xx},
\qquad
B=\partial_x,
\]
with homogeneous Dirichlet boundary conditions imposed through
the domain of \(A\).

The operators \(A\) and \(B\) do not commute as operators
with the prescribed boundary conditions. In particular, the
advection operator does not preserve the domain of the
Dirichlet Laplacian. The corresponding finite-dimensional
matrices obtained from the spatial discretization are also
noncommuting. Moreover, the operators
\(A^{-1/2}B\) and \(BA^{-1/2}\) extend to bounded operators
on \(L^2(0,1)\). Thus, Assumption~\ref{ass2} is satisfied
with \(\gamma=\frac12\).

For the spatial discretization, we employ the standard
second-order central finite difference scheme. Let
\[
x_i=ih,\qquad i=0,\ldots,N+1,
\qquad h=\frac{1}{N+1},
\]
where \(N\) denotes the number of interior grid points.
The discrete diffusion and advection operators are defined
by
\[
(A_hu)_i
=\frac{-u_{i-1}+2u_i-u_{i+1}}{h^2},
\qquad
(B_hu)_i
=\frac{u_{i+1}-u_{i-1}}{2h},
\]
for \(i=1,\ldots,N\), with \(u_0=u_{N+1}=0\).
This yields the semi-discrete system
\begin{equation}\label{eq:semidiscrete_problem2}
	u_h'(t)+A_hu_h(t)=B_hu_h(t),
	\qquad
	u_h(0)=u_{0,h},
\end{equation}
where \(u_{0,h}\) is obtained by evaluating the initial
condition at the interior grid points.

The exact solution
of the semi-discrete system is computed as
\begin{equation}\label{eq:reference_solution2}
	u_h(t)=
	\mathrm{e}^{-t(A_h-B_h)}u_{0,h}.
\end{equation}
We evaluate the corresponding matrix exponential using
MATLAB's built-in function \texttt{expm}. Consequently, the reported errors measure
only the temporal error associated with the exponential
Runge--Kutta method for the given spatial discretization;
the spatial discretization error is not included in these
comparisons.

Finally, the initial data
\[
u_0(x)=64x^3(1-x)^3
\]
is smooth and belongs to
\(\mathcal{D}((A-B)^2)\). Thus, it satisfies the corresponding
regularity requirement in Theorem~\ref{theo4}.

\begin{figure}[htbp]
	\centering
		\includegraphics[width=0.6\textwidth]{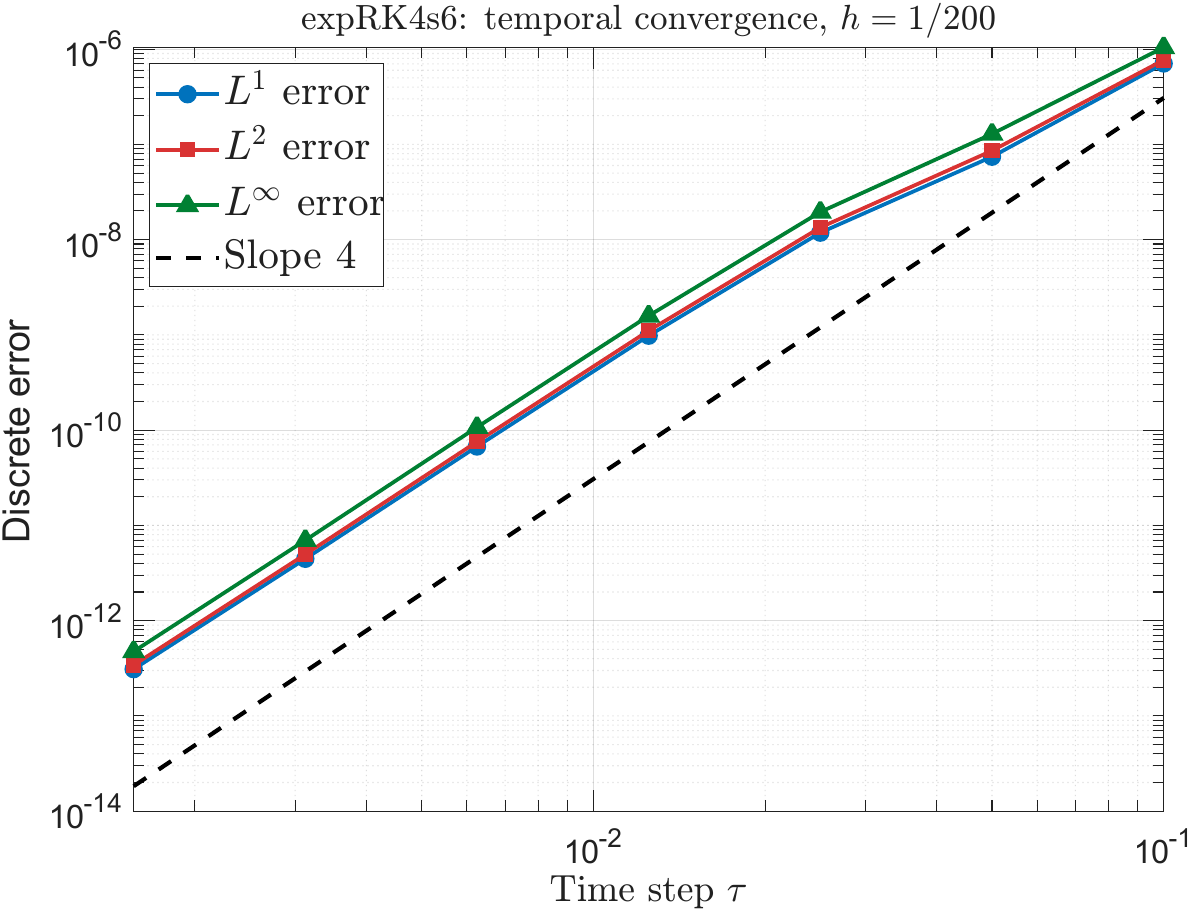}
	\caption{
		Temporal convergence of expRK4s6 in the discrete
		\(L^1(\Omega)\), \(L^2(\Omega)\), and
		\(L^\infty(\Omega)\) norms. The reference lines indicate
		fourth-order convergence.
	}
	\label{fig:convergence_fixed}
\end{figure}

We investigate the temporal convergence behaviour of the
expRK4s6 method in the discrete \(L^1(\Omega)\),
\(L^2(\Omega)\), and \(L^\infty(\Omega)\) norms. These norms
are approximated by the corresponding discrete norms
\[
\|v_h\|_{\ell_h^1}
:=h\sum_{i=1}^{N}|v_i|,
\qquad
\|v_h\|_{\ell_h^2}
:=\left(h\sum_{i=1}^{N}|v_i|^2\right)^{1/2},
\qquad
\|v_h\|_{\ell_h^\infty}
:=\max_{1\leq i\leq N}|v_i|.
\]
The time interval is fixed to \([0,T]\), with \(T=1\), and we consider the time step sizes
\[
\tau\in
\left\{
\frac{1}{10},\frac{1}{20},\frac{1}{40},
\frac{1}{80},\frac{1}{160},\frac{1}{320},\frac{1}{640}
\right\}.
\]
For a fixed spatial mesh size \(h=1/200\), Figure~\ref{fig:convergence_fixed}
presents the temporal convergence results in the three discrete norms.
The error curves exhibit consistent fourth-order convergence behaviour,
with observed convergence rates close to the reference slope of four.

These numerical observations are consistent with Theorem~\ref{theo4}.
In particular, although the theoretical estimate involves a
\(\zeta\)-dependent loss in the convergence order, \(\zeta\) can be
chosen arbitrarily small. Consequently, the theoretical convergence
rate can be made arbitrarily close to four, in agreement with the
numerical results.

To examine the influence of the spatial discretization on
the temporal convergence, we repeat the experiments for
\[
h\in\left\{\frac{1}{25},
\frac{1}{50},\frac{1}{100},
\frac{1}{200}
\right\}.
\]
For each spatial discretization, the reference solution is
computed from the corresponding semidiscrete problem using
\[
u_h(t)=\mathrm{e}^{-t(A_h-B_h)}u_{0,h}.
\]
For a given mesh size \(h\), let \(u_{h,\tau}^{n}\) denote
the numerical solution obtained with expRK4s6 at
\(t_n=n\tau\), and let \(u_h(t_n)\) denote the exact solution
of the corresponding semidiscrete system. 

\begin{figure}[h!]
	\subfigure{\includegraphics[width=0.53\textwidth]{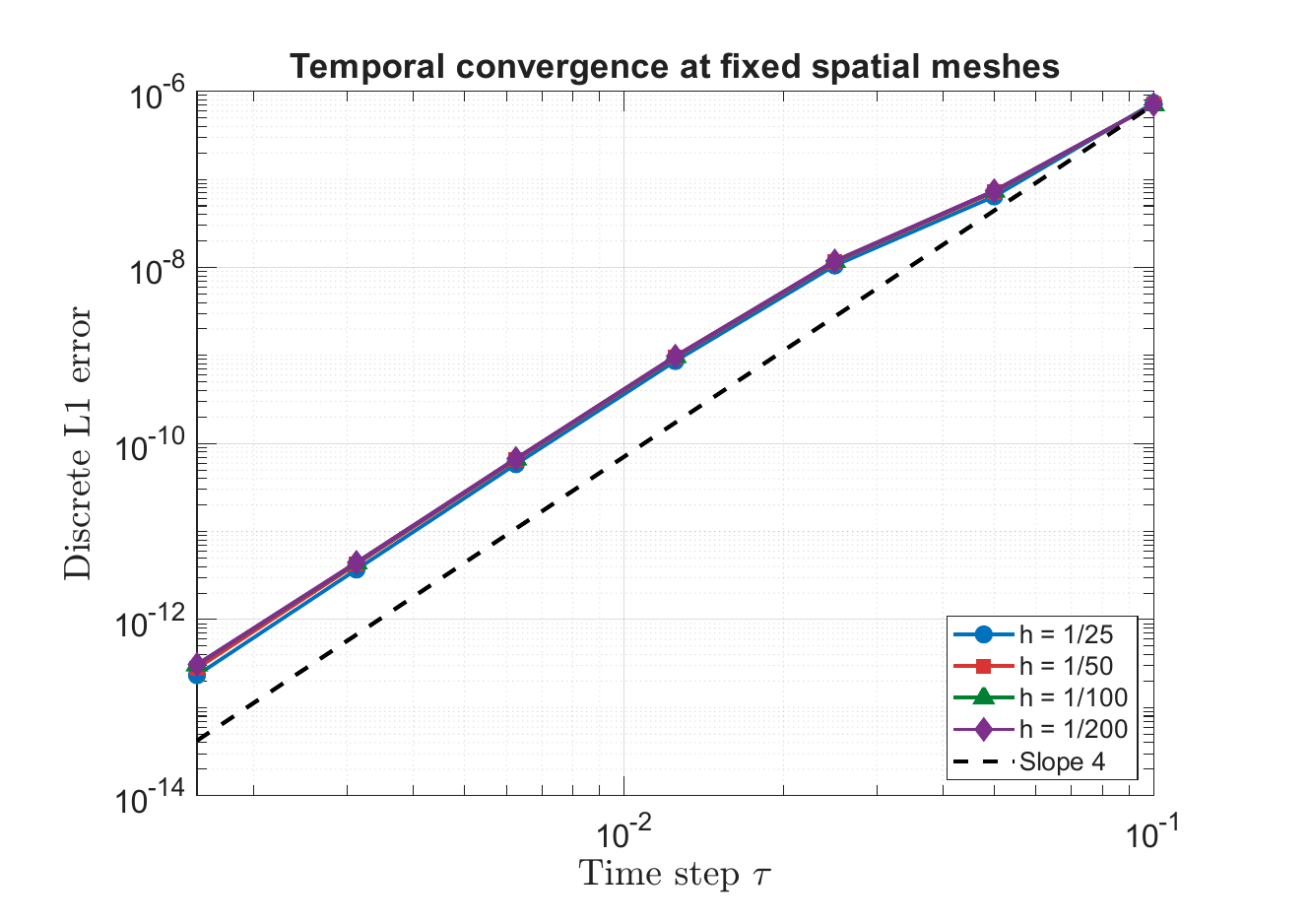}}
	\subfigure{\includegraphics[width=0.53\textwidth]{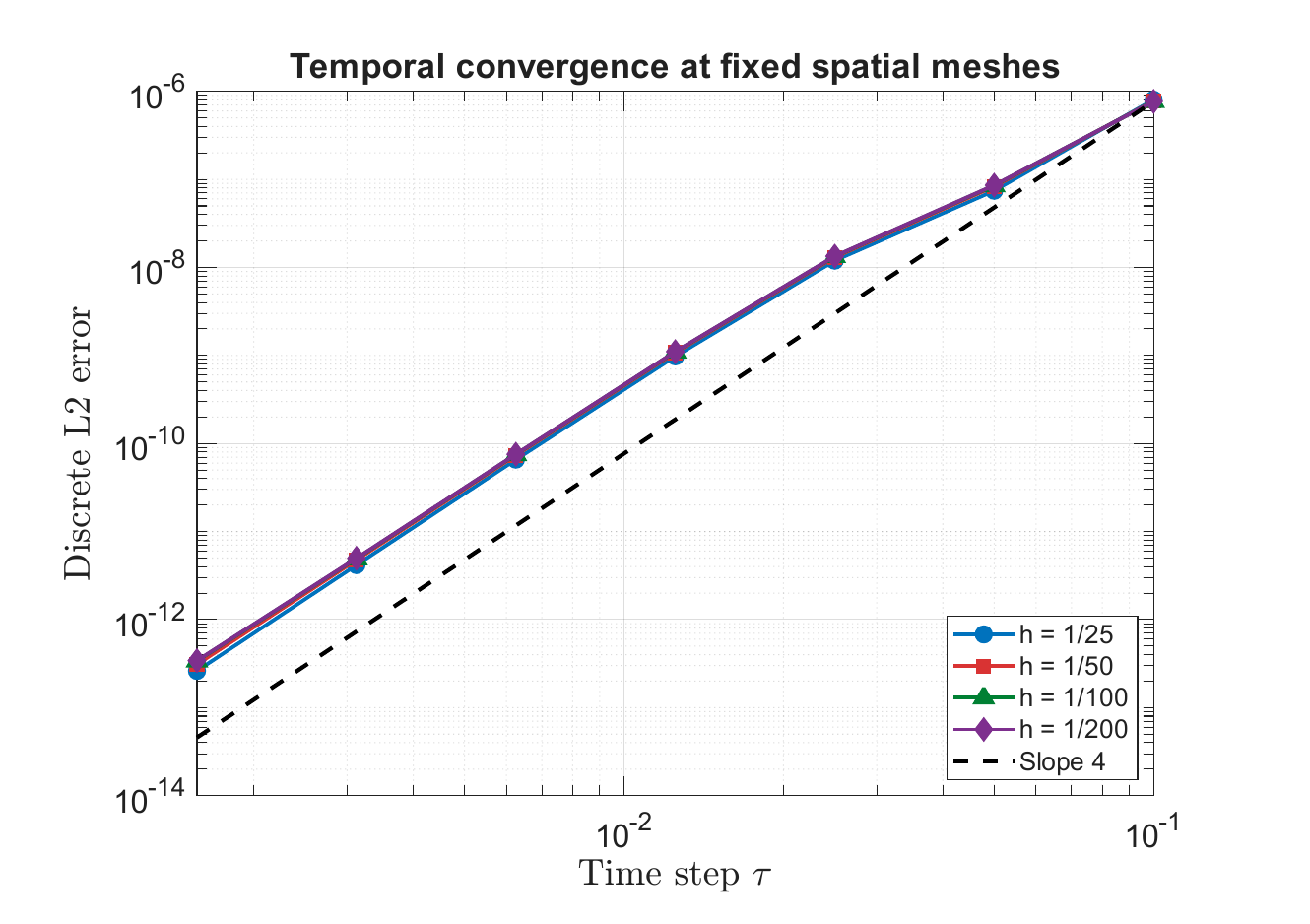}}
	\begin{center}
		\subfigure{\includegraphics[width=0.53\textwidth]{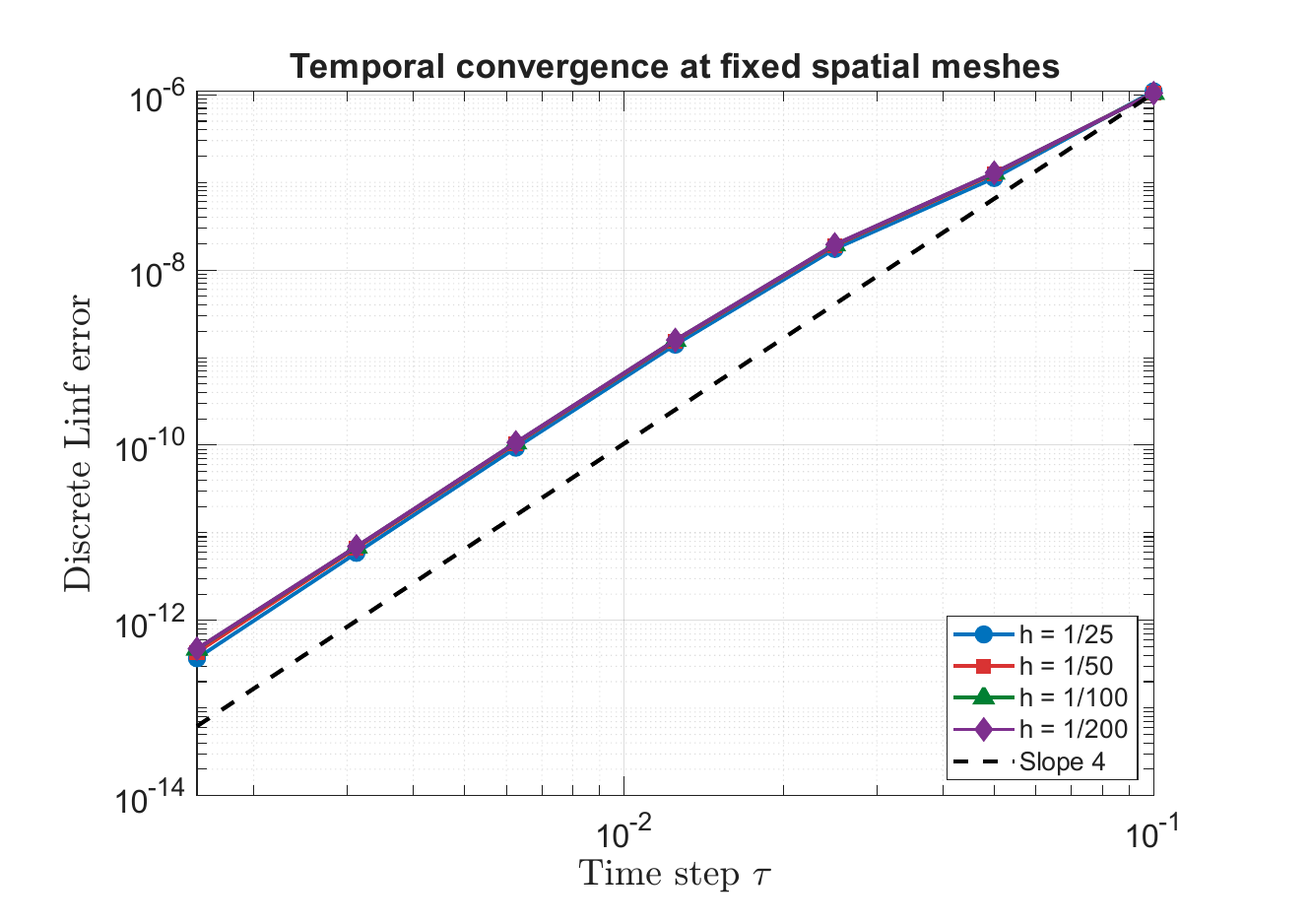}}
	\end{center}
	\caption{Temporal convergence in the discrete \(L^1(\Omega)\), \(L^2(\Omega)\),
		and \(L^\infty(\Omega)\) norms for different spatial mesh sizes \(h\).
		The dashed lines indicate the reference slope 4.}
	\label{fig:mesh_refinement}
\end{figure}

Figure~\ref{fig:mesh_refinement} presents the temporal
convergence results in the three discrete norms for the
different spatial mesh sizes. The error curves exhibit
consistent fourth-order convergence behaviour, with
observed rates close to the reference slope four. Moreover, the curves corresponding to different spatial resolutions remain parallel. This indicates that the observed temporal convergence behaviour is stable under spatial mesh refinement.

We define the
temporal error at the final time \(T=1\) by
\begin{equation}\label{eq:temporal_error}
	E_{h,p}(\tau)
	:=
	\left\|
	u_{h,\tau}^{T/\tau}-u_h(T)
	\right\|_{L^p},
	\qquad p\in\{1,2,\infty\}.
\end{equation}
Since both solutions are computed using the same spatial
discretization, this error measures exclusively the temporal
discretization error and is therefore independent of the
spatial discretization error.

To investigate the dependence of the temporal error constant
on the spatial mesh, we consider the normalized errors
\[
C_{h,p}^{\mathrm{num}}
:=
\max_{\tau\in\mathcal{T}}
\frac{E_{h,p}(\tau)}{\tau^4},
\qquad p\in\{1,2,\infty\},
\]
where \(\mathcal{T}\) denotes the set of time-step sizes
considered in the experiment and \(E_{h,p}(\tau)\) is the
temporal error defined in \eqref{eq:temporal_error}.
These quantities provide a numerical measure of the
fourth-order error constants for the different spatial
discretizations. The resulting values are reported in
Table~\ref{tab:error_constants_mesh}.

\begin{table}[H]
	\centering
	\caption{
		Normalized numerical error constants
		\(C_{h,p}^{\mathrm{num}}\) for expRK4s6 in the
		three discrete norms and for different spatial
		mesh sizes \(h\).
	}
	\label{tab:error_constants_mesh}
\begin{tabular}{lcccc}
	\toprule
	Norm & $h=1/25$ & $h=1/50$ & $h=1/100$ & $h=1/200$ \\
	\midrule
	$L^1$ &
	$0.03915$ & $0.04681$ & $0.05147$ & $0.05195$ \\
	
	$L^2$ &
	$0.04387$ & $0.05196$ & $0.05691$ & $0.05735$ \\
	
	$L^\infty$ &
	$0.06203$ & $0.07266$ & $0.07925$ & $0.07971$ \\
	\bottomrule
\end{tabular}
\end{table}
The numerical values in Table~\ref{tab:error_constants_mesh}
allow us to assess whether the normalized temporal errors
remain bounded as the spatial mesh is refined. In particular,
a stabilization of these quantities across the different
mesh sizes provides numerical evidence that the fourth-order
temporal error constant does not deteriorate under spatial
mesh refinement. Together with the convergence results in
Figure~\ref{fig:mesh_refinement}, this supports the observed
mesh-independent fourth-order temporal convergence of
expRK4s6.

\section{Conclusion}
\label{concluchap2}

In this work, we have investigated the convergence properties of the fourth-order exponential Runge--Kutta method expRK4s6 for linear parabolic problems of the form
$
u'(t)+Au(t)=Bu(t),
$
where \(A\) is an unbounded operator generating an analytic semigroup and \(B\) is a relatively bounded perturbation. Under appropriate regularity assumptions on the initial data and the operators, we have established error bounds for a general class of methods satisfying suitable order conditions, within an abstract Banach-space framework. The analysis shows that the temporal convergence rate can be made arbitrarily close to the expected fourth order.

We have also presented numerical experiments for a one-dimensional advection--diffusion problem with homogeneous Dirichlet boundary conditions. The numerical results, measured in the discrete \(L^1(\Omega)\), \(L^2(\Omega)\), and \(L^\infty(\Omega)\) norms, confirm the predicted convergence behaviour. Furthermore, experiments with different spatial mesh sizes indicate that the observed temporal convergence rates and normalized error constants remain stable under spatial mesh refinement. This provides numerical evidence for the robustness of the temporal error estimates with respect to the spatial discretization.

The present analysis contributes to the understanding of fourth-order exponential Runge--Kutta methods for problems involving non-commuting unbounded operators. Further research may investigate higher-order methods, weaker regularity assumptions on the initial data, and extensions to more general parabolic problems and spatial discretizations.

\subsection*{Funding}
Thi Tam Dang has been supported by the Jane and Aatos Erkko Foundation.
Trung Hau Hoang gratefully acknowledges support from Charles
University through the research programmes PRIMUS/25/SCI/025
and UNCE/24/SCI/005.

\subsection*{Conflict of interest}
The authors declare that there are no conflicts of interest.

\subsection*{Data availability}
No data was used for the research described in the article.

\end{document}